\documentclass[12pt]{article}
\title{A numerical study of the pull-in instability
 in some free boundary models for MEMS}

\author{Gilberto Flores,\\ 
Instituto de Investigaciones en Matem\'aticas 
Aplicadas y en Sistemas, \\ and FENOMEC\\ 
Universidad Nacional Aut\'onoma de 
M\'exico, \\Apdo. Postal 20-126, \\
01000 M\'exico, D.F., MEXICO  \\
(gfg@mym.iimas.unam.mx)
\and Noel F. Smyth,\\
School of Mathematics, \\ University of Edinburgh, 
\\ Edinburgh, Scotland, U.K., EH9 3FD.\\
(N.Smyth@ed.ac.uk)}
\date{}
\usepackage[dvips]{graphicx}

\def\complejo{\hbox{{\sf C}\kern-.45em\lower-.45ex\hbox{{\tiny ]}}\enskip}}
\def\real{\hbox{{\sf I}\kern-.1em{\sf R}}}
\def\natural{\hbox{{\sf I}\kern-.1em{\sf N}}}
\def\realc{\hbox{{\tiny I}\kern-.1em{\tiny R}}}
\def\entero{\hbox{{\sf Z}\kern-.40em{\sf Z}}}

\begin{document}

\baselineskip=.6truecm

\maketitle







\begin{abstract}

In this work we numerically compute the bifurcation curve of stationary
solutions for the free
boundary problem for MEMS in one space dimension. It has a single turning
point, as in the case of the small aspect ratio limit. We also find a 
threshold for the existence of global-in-time 
solutions of the evolution equation given by either a heat 
or a damped wave equation.  This threshold is what we term the dynamical 
pull-in value: it separates the stable operation regime from the touchdown 
regime.  The numerical calculations show that the dynamical threshold 
values for the heat equation coincide with the static values.
For the damped wave equation the dynamical threshold values are smaller 
than the static values. 
This result is in agreement with the observations reported for a 
mass-spring system studied in the engineering literature.  
In the case of the damped wave equation, we also show that 
the aspect ratio of the device is more important than the inertia in the 
determination of the pull-in value. 

\end{abstract}

\newpage

{\bf 1. Introduction}

The operation of many micro electromechanical systems (MEMS) relies upon the
action of electrostatic forces.  Many such devices, including pumps, switches or
valves, can be modelled by electrostatically deflected elastic membranes.
Typically a MEMS device consists of an elastic membrane held at a constant 
voltage and suspended above a rigid ground plate placed in series with a 
fixed voltage source.  The voltage difference causes a deflection of the 
membrane, which in turn generates an electric field in the region between 
the plate and the membrane.  Mathematically, this is then a free boundary 
problem.  The electric potential is defined in a region which depends on the 
membrane deflection, while the elastic deformation is forced by the trace of 
the electric field on the membrane.

An important nonlinear phenomenon in electrostatically deflected membranes
is the so-called ``pull-in'' instability.  For moderate voltages the system
is in the stable operation regime: the membrane approaches a steady state
and remains separate from the ground plate.  When the voltage is increased
beyond a critical value, there is no longer an equilibrium configuration
of the membrane.  As a result, the membrane collapses onto the ground
plate.  This phenomenon is also known as ``touchdown''.  The critical value
of the voltage required for touchdown to occur is termed the 
{\it pull-in value}.  
The determination of the pull-in value is important for the design 
and manufacture of MEMS devices, particularly as touchdown is a desirable 
property in devices such as microvalves.  For instance, Desai {\em et al} 
\cite{desai} give a description of microvalves used in microfluidic chips. 
However, for most devices, it is desirable to achieve the stable operation 
regime with no touchdown. The {\it pull-in distance} is the critical
distance between the ground plate and the elastic membrane beyond which
pull-in occurs.

The issue of the static and dynamical pull-in instabilities has been 
addressed by 
the engineering community in the context of a model in which the moving 
structure 
is a plate attached to a spring with damping.  The elastic properties of 
the moving 
plate are described by the restoring force of the spring, which is assumed to 
be 
given by Hooke's law. The voltage applied to the moving plate results in an
electrostatic force acting on the spring-mass system,  see Rocha {\em et al} 
\cite{rocha} and Zhang {\em et al} \cite{zhang} for details.  
The governing equation
for the displacement of the moving mass is
\begin{equation}
m {d^{2}x \over dt^{2}} + b {dx\over dt} + kx = {\lambda \over (d_0 -x)^{2}},
\label{e:spring}
\end{equation}
where $d_0$ is the initial gap between the plates, $\lambda= \epsilon_{0} A V^{2}/2$, 
$A$ is the area of the moving plate, $V$ is the voltage applied to it and $\epsilon_{0}$ is 
the permitivity of the free space between the plates. The right hand side is the
Coulomb force.

Zhang {\em et al} \cite{zhang} described the dynamical pull-in as the 
collapse of the moving structure towards the substrate due to the combined
action of kinetic and potential energies.  They also stated that, in general,
{\bf dynamical pull-in requires a lower voltage to be triggered compared to
the static pull-in threshold}.  One of the findings in Rocha {\em et al.\ } 
\cite{rocha} is the fact that for an overdamped device, the dynamics
in the  touchdown regime
has three distinguished regions characterized by different time scales:
in the first region the structure moves fast until it gets near the static
pull-in distance, at which point there is a metastable region in which the 
motion is very slow, and finally a third region in which collapse takes place 
on a fast time scale.  

In the present work we study the equations obtained when the deflection
of the elastic membrane is governed by a forced, damped wave equation. 
We  also study the forced
heat equation which corresponds to setting the inertia equal to zero.
For simplicity, we assume that the motion starts from rest.
The numerical results indicate that for the forced heat equation, the dynamical 
pull-in value coincides with the static value. 
This result is supported by the fact that the
membrane profiles decrease monotonically in time and approach a steady
state in the stable operation regime, which suggests that there is a 
maximum principle and that the stationary solutions act as a barrier to
prevent touchdown. This is exactly the situation in the limit of
vanishing aspect ratio. In contrast, the dynamical pull-in value for the 
damped wave equation 
is smaller than the static value, in agreement with
the observation in \cite{zhang} for equation (\ref{e:spring}). 
We also obtain the different time scales in the dynamics of
touchdown as reported in \cite{rocha} for equation (\ref{e:spring}).  
Our results then indicate that the difference between the dynamical and
static pull-in values is due to the inertial forces. 
On the other hand, our calculations show that the 
aspect ratio is more important than the
inertia in the determination of the dynamical pull-in value.

Here we study the following free-boundary problem.
Let $u$ denote the membrane deformation.  In terms of dimensionless 
variables, the electric potential $\psi$ is defined in the region 
\begin{equation} 
\Omega(u) = \{(x,z) \in (-1,1) \times (-1, \infty) :\, -1<z<u(x) \} .
\label{e:region}
\end{equation}
The electric potential itself is the solution of the elliptic equation
\begin{equation}
   \epsilon^{2} \, \psi_{xx} + \psi_{zz}=0 
\label{e:poteqn}
\end{equation}
together with the boundary conditions
\begin{equation}
\psi(x,-1)=0, \, \, \psi(x, u(x))=1 \, \, {\rm for} \, \, x \in (-1,1),
\, \, \psi(\pm 1, z)= 1+z \, \, {\rm for} \, \, z \in (-1,0).
\label{e:bc}
\end{equation}
The membrane deformation $u$ itself is the solution of
\begin{equation}
\gamma \,  u_{tt} + u_{t} - u_{xx}= -\lambda \, [ \epsilon^{2} \,
| \psi_{x} (x,u(x))|^{2} + |\psi_{z} (x,u(x))|^{2}] .
\label{e:mem}
\end{equation}
For simplicity, we assume that the motion starts from the rest position.
In these equations, the control parameter $\lambda$ is proportional to the 
square of the applied voltage, $\epsilon$ is the ratio of the gap size to 
the device length and $\gamma$ is the ratio of inertial to damping forces.
For a derivation of these equations see Pelesko and Bernstein \cite{pelesko}.

In the formulation above there are other effects which have not been included.
One is the effect of the electric field at the edge of the membrane, known as 
{\it fringing fields}.   In addition, the elastic energy in the present 
model does not include the curvature of the membrane.  Pelesko and Driscoll 
\cite{pelesko2} gave a derivation of the governing equation when the fringing 
field is taken into account.  The boundary value problem for the electric potential 
(\ref{e:poteqn})--(\ref{e:bc}) is then solved for $\epsilon =0$ with a boundary
layer correction around the edge of the membrane.  Brubaker and Pelesko \cite{brubaker} 
studied the case in which the elastic energy includes the curvature 
of the membrane. The electric potential is obtained for $\epsilon=0$.

The small aspect ratio limiting case corresponding to $\epsilon=0$ has been 
studied extensively.  In this case, the boundary value problem (\ref{e:poteqn}) and 
(\ref{e:bc}) for the electric potential can be solved explicitly to give
\begin{equation}
\psi (x,z)= \frac{1+z}{1+u(x)}, \quad (x,z) \in \Omega(u).
\label{e:limitpot}
\end{equation}
Equation (\ref{e:mem}) for the elastic deformation then reduces to a nonlinear
wave equation, termed the small aspect ratio model
\begin{equation}
\gamma u_{tt} + u_{t} - u_{xx} = -\frac{\lambda}{(1+u)^{2}}.
\label{e:aspect}
\end{equation}
The further limiting case with $\gamma=0$ is a nonlinear heat equation for 
which the dynamical pull-in value coincides with the critical value 
$\lambda^{*}$ for the existence of stationary solutions of (\ref{e:aspect}).  Indeed, 
there are two, one or zero stationary solutions of (\ref{e:aspect}) according to 
whether $\lambda < \lambda^{*}$, $\lambda = \lambda^{*}$ or 
$\lambda > \lambda^{*}$.  Moreover, solutions of the nonlinear
heat equation corresponding to $\lambda \leq \lambda^{*}$ converge 
to a steady state, while solutions corresponding to $\lambda > \lambda^{*}$ 
quench in finite time.  The proof of this behaviour relies on the maximum 
principle, see Flores {\em et al} \cite{flores2}.

The same behaviour is obtained when the effect of fringing fields is taken 
into account.  According to Pelesko and Driscoll \cite{pelesko2},
equation (\ref{e:aspect}) is modified as follows. 
The numerator on the right hand side becomes $\displaystyle -\lambda
(1+ \epsilon^{2} u_{x}^{2})$.  For stationary solutions, Lindsay and
Ward \cite{lindsay} have established that the pull-in value 
$\lambda^{*}(\epsilon)$ admits an asymptotic expansion in 
powers of $\epsilon^{2}$ and obtained the leading order 
term, which in the one-dimensional case corresponds to the critical value
$\lambda^{*}$ mentioned in the previous paragraph. 
Wei and Ye \cite{wei} have described the structure of the stationary
solutions for this problem.  There is a critical value of $\lambda$
such that there are at least two solutions, one or none according to whether
$\lambda$ is smaller, equal to, or larger than this critical value.
Liu and Wang \cite{liu} verified that for the corresponding heat equation, 
the dynamical critical parameter coincides with the static critical value.  
The stationary solution thus acts as a barrier and prevents touchdown.
The rule is that the static and dynamical pull-in values coincide whenever 
there is a maximum principle.

On the other hand, the numerical evidence for the case $\gamma > 0$
indicates that for the damped wave equation (\ref{e:aspect}) 
there is a threshold, which we denote by $\lambda^{*}_{w}$,
that separates the stable operation regime from the touchdown regime. 
This means that solutions of (\ref{e:aspect}) converge to a steady state for 
$\lambda<\lambda^{*}_{w}$, while for $\lambda>\lambda^{*}_{w}$ solutions
quench in finite time.  This critical value of $\lambda$ is what we call the 
{\it dynamical pull-in value}.  Moreover, $\lambda^{*}_{w} < \lambda^{*}$, see 
Flores \cite{flores}.  Similar numerical results concerning the dynamical 
threshold were obtained for conservative wave equations with a singular forcing
term in one dimension by Chang and Levine \cite{chang} and in higher dimensions 
by Smith \cite{smith}.  In the same context, Kavallaris {\em et al} \cite{kavallaris} 
numerically found the existence of a dynamical threshold, smaller than the 
static value, in a one dimensional, non-local 
version of the equation considered in \cite{chang} for which the MEMS device 
is connected in series with a capacitor.

The experimental investigation of Siddique {\em et al} \cite{siddique} 
points in the same direction.  They set up an array of two plates, one 
fixed, the other with a laser cut hole where a soap film was applied. 
The plates were separated by a distance $d$.  The critical voltage was 
computed for different values of $d$.  An empirical relation was then used to 
determine the critical value of $\lambda$.  These values were compared with either 
upper and lower bounds or with numerically computed values of $\lambda^{*}$ for 
elliptical or rectangular domains.  Good agreement was found for small values of $d$.  
It was found that the experimental values were smaller than the numerically calculated values.  
The interpretation of this is that the experimental values correspond to 
the dynamical critical value of $\lambda$.  In Siddique {\em et al.\ } \cite{siddique} a 
question is raised so as to identify the most important effect which accounts for 
the difference between the theoretical and the experimental results.  The numerical results 
of the present work indicate that the aspect ratio of the device is more important than the 
inertial effects.  Another part of the explanation is that the static and 
dynamical pull-in values are different.  

The static free boundary problem and the associated semilinear parabolic 
equation in one space dimension governing it have been analyzed by 
Lauren\c{c}ot {\em et al} \cite{lauren} and by Escher {\em et al} 
\cite{escher}, respectively.  In the first work the existence of stationary solutions for 
small values of $\lambda$ was established, as well as the non-existence for large 
values of this control parameter. The local well-possedness of the parabolic
problem was proved in \cite{escher}.  It was also established that for small values of 
$\lambda$ the solution exists for all times and converges to a steady state as $t \to \infty$.  
It was also proved that for large values of $\lambda$ global existence does not hold in 
the sense that $u$ reaches the value $-1$ in finite time, that is, $\displaystyle u$
{\it quenches} in finite time.  To the best of our knowledge these are 
the only rigorous results to date for the free boundary problem.
As mentioned in Lauren\c{c}ot {\em et al} \cite{lauren}, no further 
information is available on the structure of the set of values of 
$\lambda$ for which there is a classical stationary solution of the 
free boundary problem.  It is believed that this set is an interval.
In the present work, by computing the bifurcation curve we provide numerical 
evidence that this is indeed the case.  The shape of the bifurcation curve 
for the steady states is qualitatively similar to the corresponding curve 
for the small aspect ratio limit corresponding to $\epsilon=0$, which 
suggests the existence of a critical value $\lambda^{*}(\epsilon)$ for a 
steady state to exist.  The numerical results also indicate that 
$\lambda^{*}(\epsilon)\to \lambda^{*}$ as $\epsilon \to 0^{+}$.

We also provide numerical evidence which shows that this static critical 
value $\lambda^{*}(\epsilon)$ coincides with the dynamical pull-in value
for the nonlinear heat equation.  In contrast, for the damped wave equation it
does not control the dynamics since the dynamic pull-in value is smaller
than the static critical value, even in the limiting case $\epsilon=0$.
Therefore, the difference between the dynamic and static critical values
is due to the inertial forces.  We also find that the aspect ratio $\epsilon$ is 
more important than the inertia coefficient $\gamma$ in the determination of the 
dynamical pull-in value.

\section {Stationary solutions}

As discussed in the previous section, the equation for the electric 
potential $\displaystyle \psi$ is
\begin{equation}
\epsilon^{2} \, \psi_{xx} + \psi_{zz}=0 
\label{e:potpsi}
\end{equation}
in the region $\Omega(u) = \{(x,z) \in (-1,1) \times (-1, \infty) :
\, -1<z<u(x) \}$, together with the boundary conditions 
\begin{equation}
\psi(x,-1)=0, \, \, \psi(x, u(x))=1 \, \, {\rm for} \, \, x \in (-1,1),
\quad \psi(\pm 1, z)= 1+z \, \, {\rm for} \, \, z \in (-1,0).
\label{e:potbc}
\end{equation}
The elastic deformation $\displaystyle u$ is the solution of
\begin{equation}
u_{xx}= \lambda \, [ \epsilon^{2} \,
| \psi_{x} (x,u(x))|^{2} + |\psi_{z} (x,u(x))|^{2}] 
\label{e:elasticu}
\end{equation}
with the boundary condition $u(\pm 1)=0$.


Following Lauren\c{c}ot {\em et al} \cite{lauren}, we map the 
domain $\Omega(u)$ onto the rectangle
\begin{equation}
\Omega = (-1,1) \times (0,1)
\label{e:newdom}
\end{equation}
by means of the transformation
\begin{equation}
T_{u} (x,z) = \left(x, \frac{1+z}{1+u(x)}\right),
\label{e:trans}
\end{equation}
which has the inverse
\begin{equation}
T_{u}^{-1}(x, \eta)= (x, [1+u(x)]\eta -1).
\label{e:intran}
\end{equation}
In terms of the new independent variables $(x, \eta)$, the electric potential
is denoted by $\phi$: $\phi= \psi \circ T_{u}^{-1}$.  The potential equation
(\ref{e:potpsi}) then becomes
\begin{equation}
{\cal L}_u(\phi) =0 \quad {\rm in} \quad \Omega , \qquad \phi(x,\eta) = \eta \quad 
{\rm on} \quad \partial \Omega  ,
\label{e:transpot}
\end{equation}
where $\displaystyle {\cal L}_u$ is the elliptic operator defined by
\begin{equation}
{\cal L}_{u}(\phi) = \epsilon^{2} \phi_{xx} -2 \epsilon^{2} \eta
\frac{u_{x}}{1+u(x)}\phi_{x \eta} + 
\frac{1+\epsilon^{2}\eta^{2} u_{x}^{2}}{ 
[1+u(x)]^{2}} \phi_{\eta \eta} +
 \epsilon^{2} \eta \left[ 2 \left(\frac{u_x}{ 1+u(x)}\right)^{2}
-\frac{u_{xx}}{ 1+u(x)} \right] \phi_{\eta}.
\label{e:oper}
\end{equation}
Equation (\ref{e:elasticu}) for the elastic deformation $u$ becomes 
\begin{equation}
u_{xx} = \lambda \left[\frac{1 + \epsilon^{2} u_{x}^{2}}{(1+u(x))^{2}}\right]
|\phi_{\eta} (x,1)|^{2} 
\label{e:ueqntrans}
\end{equation}
in $(-1,1)$, with the boundary condition $u(\pm 1)=0$.



\begin{figure}
\centering
\includegraphics[scale=0.25,angle=270]{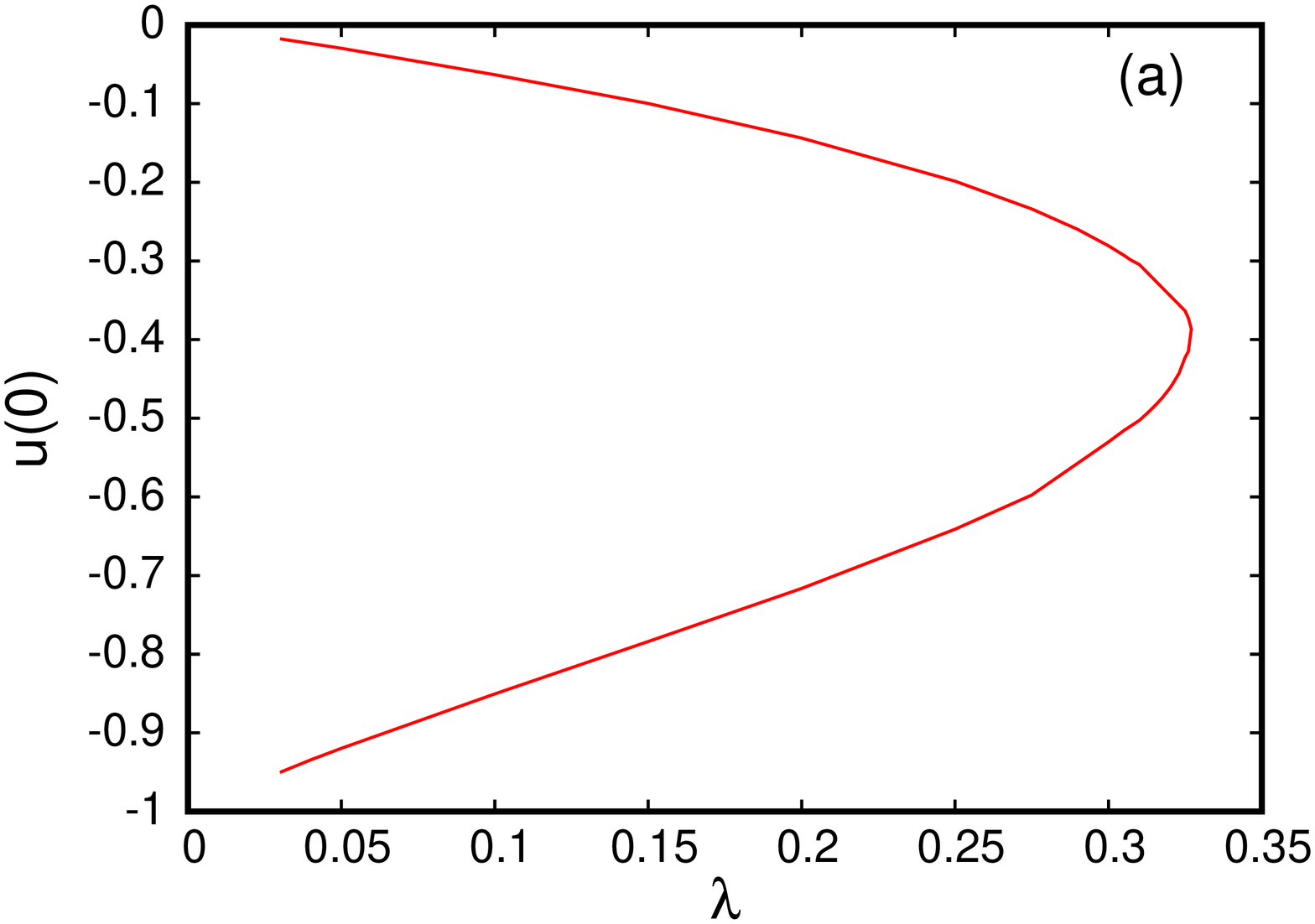}
\includegraphics[scale=0.4,angle=270]{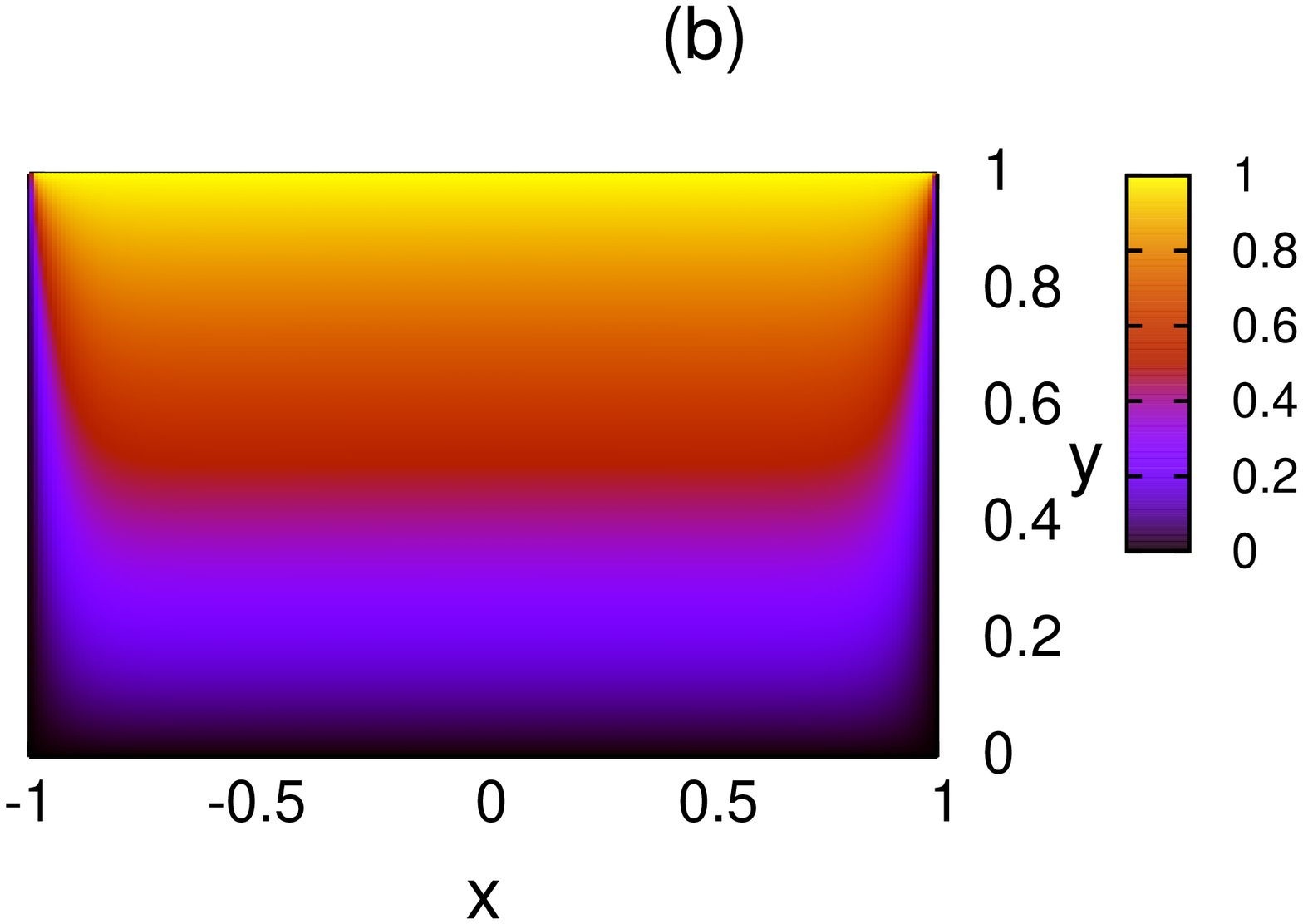}
\caption{(Color online)  $\epsilon = 0.2$.  (a) Birfurcation diagram
from numerical solution of steady equations (\ref{e:transpot}) and
(\ref{e:ueqntrans}). (b) contour plot of $\phi$ for $\lambda =0.32$.}
\label{f:bifurc}
\end{figure}

The transformed potential and elastic equations (\ref{e:transpot}) and 
(\ref{e:ueqntrans}) were solved numerically using centred finite differences
for the derivatives, so that the errors are $O(\Delta x^{2},\Delta \eta^{2})$.  
The potential equation (\ref{e:transpot}) then becomes a linear system in $\phi$ 
which was solved using Jacobi iteration.  The elastic equation (\ref{e:ueqntrans}) 
is a nonlinear two point boundary value problem and was solved using a shooting method.  
The potential equation (\ref{e:transpot}) and the elastic equation (\ref{e:ueqntrans}) 
form a coupled system due to $u$ appearing in the elliptic operator (\ref{e:oper}).  
A Picard iteration was then used to solve this coupled system.  A starting
guess for $\phi_{\eta}(x,1)$ was assumed and then the elliptic equation (\ref{e:transpot}) was
solved to find $\phi$ and so $\phi_{\eta}$ at $\eta=1$.  The deformation equation (\ref{e:ueqntrans})
was then solved for $u$ using this $\phi_{\eta}(x,1)$.  With this updated 
$u(x)$ the elliptic equation (\ref{e:transpot}) was again solved and the process iterated until 
convergence.  The numerical results show the existence of a critical value of $\lambda$, 
denoted by $\lambda_{s}^{*}(\epsilon)$, such that there are two, one or zero stationary solutions 
according to whether $\lambda$ is below, equal to or above the critical value $\lambda_{s}^{*}(\epsilon)$.  
A low initial guess for $u'(-1)$, between $0$ and $-1.5$, resulted in the numerical solution
for $u$ converging to the upper branch of solutions and a high initial guess for $u'(-1)$, between $-1.5$ 
and $-3$, resulted in convergence to the lower branch.  For $\epsilon=0$, it is known that 
$\lambda_{s}^{*} = 0.350004\ldots$ \cite{flores2}.  The numerical scheme was tested by finding 
$\lambda^{*}_{s}$ in the limit $\epsilon \to 0$ with $\Delta x = \Delta \eta = 5\times 10^{-3}$.  
For $\epsilon = 0.0001$ it was found that $\lambda_{s}^{*} = 0.350000$, which agrees with the value 
for $\epsilon=0$ to five decimal places, which is the accuracy for the critical $\lambda$ which will 
be used in this work.  The bifurcation curve for $\epsilon=0.2$ is shown in Figure \ref{f:bifurc}(a). 
The bifurcation parameter chosen was the value of $u$ at $x=0$.  
Figure \ref{f:bifurc}(b) shows a contour plot of the electric potential
$\phi$.  Due to $\epsilon$ being small, over a large part of the domain the electric potential 
for the free boundary problem is close to the potential for the small 
aspect ratio limit (\ref{e:limitpot}), 
which in the transformed variables is $\phi_0(x,\eta)=\eta$.

\section{Dynamical solutions}

The dynamical behaviour of the membrane was also investigated, as discussed above for
the small aspect ratio model (\ref{e:aspect}).  To investigate the dynamical behaviour
of the membrane, the forced heat equation 
\begin{equation}
u_{t} - u_{xx} = -\lambda \left[\frac{1 + \epsilon^{2} u_{x}^{2}}{(1+u(x))^{2}}\right]
|\phi_{\eta} (x,1,t)|^{2} 
\label{e:ueqntheat}
\end{equation}
and the forced, damped wave equation
\begin{equation}
\gamma u_{tt} + u_{t} - u_{xx} = -\lambda 
\left[\frac{1 + \epsilon^{2} u_{x}^{2}}{(1+u(x))^{2}}\right]
|\phi_{\eta} (x,1,t)|^{2} 
\label{e:ueqntwave}
\end{equation}
were solved for the membrane displacement $u$.  As mentioned in the 
Introduction, we assume that the motion starts from rest. This means that the 
initial condition for the heat equation (\ref{e:ueqntheat}) is 
$u(x,0)=0$, while for the damped wave equation (\ref{e:ueqntwave}) we take 
$u(x,0)=0$ and $u_{t}(x,0)=0$. 

The forced heat equation (\ref{e:ueqntheat}) was solved using centred differences in 
space $x$ and an Euler scheme in time $t$, resulting in an explicit scheme with error 
$O(\Delta t)$ in time and error $O(\Delta x^{2},\Delta \eta^{2})$ in space, the same
spatial accuracy as the numerical scheme used to solve the potential equation 
(\ref{e:transpot}) and which was discussed in the previous section.  The same Picard
iteration as discussed in the previous section was used to find $\phi_{\eta}(x,1)$ in the 
deformation equation (\ref{e:ueqntrans}).  Except for the first time step, the value
of $\phi_{\eta}(x,1)$ at the previous time step was used as the initial guess for 
the iteration.  The potential equation (\ref{e:transpot}) was again solved using 
Jacobi iteration.  The solution for $\phi$ at the previous time step was used as the 
initial guess, which resulted in fast convergence.  The accuracy of the heat equation was 
again tested by finding the critical $\lambda$ in the limit $\epsilon \to 0$ as in this limit the 
heat equation must give the known value $\lambda^{*} = 0.350004\ldots$ \cite{flores2}.  For
$\epsilon = 0.0001$, $\Delta t = 1 \times 10^{-5}$ and $\Delta x = \Delta \eta = 5 \times 10^{-3}$ 
the critical value $0.350006$ was found, which agrees with $\lambda_{s}^{*}$ to five
decimal places.  Note that the electric potential now depends on time due to the time dependence
of the coefficients of the elliptic operator $\displaystyle {\cal L}_u$ defined in (\ref{e:oper}).

\begin{table}
\begin{center}
\begin{tabular}{|c|c|c|}  \hline
 $\epsilon$ & $\lambda^{*}_{s}$ static equation & 
$\lambda^{*}_{h}$ heat equation 
\\ \hline
0.01        &               0.34997          &              0.34996 \\ \hline
0.1         &               0.34536         &               0.34535 \\ \hline
0.2         &               0.32738         &               0.32736 \\ \hline
0.3         &               0.29356         &               0.29353 \\ \hline
\end{tabular}
\end{center}
\caption{Critical values $\lambda^{*}_{s}$ for stationary solution for steady
equations (\ref{e:transpot}) and (\ref{e:ueqntrans}) (second column) and
$\lambda^{*}_{h}$ obtained from the potential equation (\ref{e:transpot}) and 
forced heat equation (\ref{e:ueqntheat}) (third column).}
\label{t:critheat}
\end{table}

The forced, damped wave equation (\ref{e:ueqntwave}) was solved using centred differences 
in space $x$ and time $t$, again resulting in an explicit scheme with error $O(\Delta t^{2})$
in time and $O(\Delta x^{2},\Delta \eta^{2})$ in space, again the same spatial accuracy as
the scheme used to solve the potential equation (\ref{e:transpot}).  The same Picard
iteration as for the stationary solutions of the previous section and the solution of the
heat equation was used to find $\phi_{\eta}(x,1)$ from the elastic equation (\ref{e:ueqntrans}) 
with the iteration started with the value of $\phi_{\eta}(x,1)$ at the previous time step, as for 
the heat equation.  As for the heat equation, the potential equation (\ref{e:transpot}) was
solved using Jacobi iteration, as using the solution at the previous time step as the initial
guess resulted in fast convergence.  The scheme was tested by decreasing the space and time steps 
until the critical values of $\lambda$ did not change to five decimal places.  It was found
that $\Delta t = 2 \times 10^{-3}$ and $\Delta x = \Delta \eta = 5\times 10^{-3}$ were sufficient
for this.

\begin{figure}
\centering
\includegraphics[scale=0.4,angle=270]{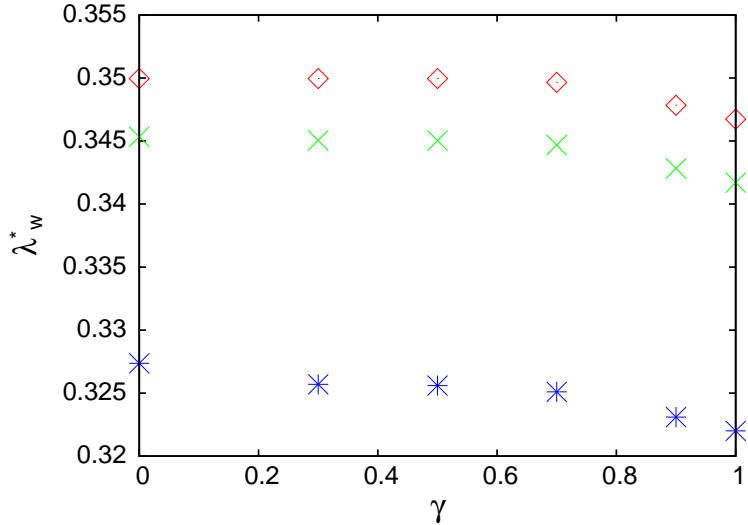}
\caption{(Color online) Plot of the critical values of $\lambda$ as
a function of $\gamma$ for $\epsilon=0.01$: (red) diamonds, $\epsilon=0.1$:
(green) cross, $\epsilon=0.2$: (blue) star.}   
\label{f:critwaveplot}
\end{figure}

\begin{figure}
\centering
\includegraphics[scale=0.4,angle=270]{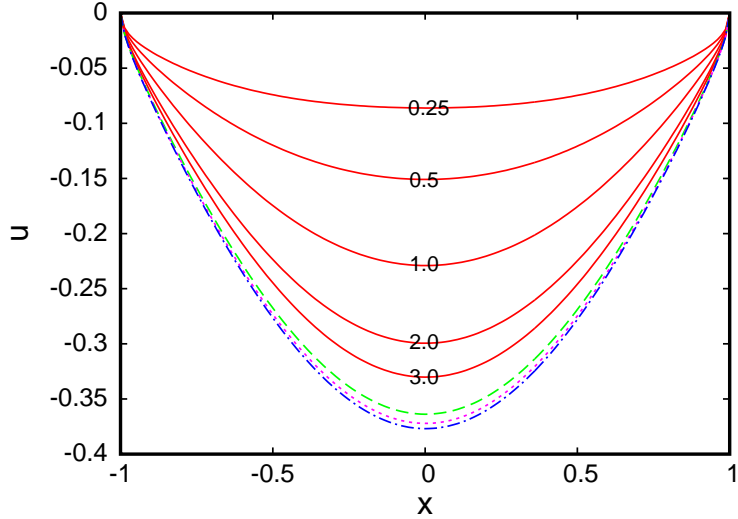}
\caption{(Color online)  Solution for $u(x,t)$ for the potential equation
(\ref{e:transpot}) and the heat equation (\ref{e:ueqntheat}) for $\epsilon = 0.2$ and
$\lambda=0.327$.  The times $t$ for each solution are on the solution curve, except for $t=6$:
green (long dash) curve, $t=8$: pink (short dash) curve and $t=10$: blue (dot-dash) curve.}
\label{f:heatevolve}
\end{figure}

The dependence of the critical value $\lambda^{*}_{s}(\epsilon)$ for 
a steady solution $u$ to exist is further illustrated in Table
\ref{t:critheat}, with the dynamic critical values illustrated in Table \ref{t:critheat} and 
Figure \ref{f:critwaveplot}.  The table and figure show the critical values  
$\lambda^{*}_{s}$, $\lambda^{*}_{h}$, and $\lambda^{*}_{w}$ as found from the steady equations 
(\ref{e:transpot}) and (\ref{e:ueqntrans}), the potential equation (\ref{e:transpot}) and the 
forced heat equation (\ref{e:ueqntheat}) for $u$ and the potential equation (\ref{e:transpot})
and the forced, damped wave equation (\ref{e:ueqntwave}) for $u$, respectively.  
As discussed above, the dynamical critical value $\lambda^{*}_{h}$ as 
determined from the forced heat equation for $u$ is slightly lower than the static value. 
However, the difference is so small and the monotonic in time behaviour of the membrane
profiles $u$ make us believe that the two critical values are equal.  The monotonic approach of 
$u$ to the steady state when the elastic deformation is governed by the forced heat equation is 
illustrated in Figure \ref{f:heatevolve}.  Note that by $t=10$ the solution has reached the steady 
state.  Note that in Figure \ref{f:critwaveplot} the values $\lambda^{*}_{h}$ have been plotted as 
the points with $\gamma = 0$.  To summarize the results, we have that for $\epsilon>0$, $\displaystyle
\lambda^{*}_{w}(\epsilon) < \lambda^{*}_{h}(\epsilon) = \lambda^{*}_{s}(\epsilon)$.  In the case 
$\epsilon=0$, it is known that $\displaystyle \lambda^{*}_{h} = \lambda^{*}$, while the 
results of Flores \cite{flores} indicate that $\lambda^{*}_{w} < \lambda^{*}$.  


The behaviour of the pull-in value is more involved when the displacement $u$ 
is given by the damped, forced wave equation (\ref{e:ueqntwave}), as can be seen
on comparing the critical values in Table \ref{t:critheat} for the heat equation 
and Figure \ref{f:critwaveplot} for the damped wave equation, again noting that the 
values $\lambda^{*}_{s}$ have been plotted as the points for $\gamma = 0$ in Figure
\ref{f:critwaveplot}.  For low values of the inertia $\gamma$ the critical value 
$\lambda^{*}_{w}$ is little changed from $\lambda^{*}_{h}$.  This is to be expected as 
the damping $u_{t}$ dominates the inertia term $\gamma u_{tt}$ in the forced, damped 
wave equation (\ref{e:ueqntwave}) for small inertia coefficient $\gamma$.
There is little change in the critical value $\lambda^{*}_{w}$ for $\gamma$
up to $0.5$.  Increasing the inertia $\gamma$ to $0.7$ results in a
significant change in $\lambda^{*}_{w}$ over $\lambda^{*}_{h}$, with the former 
value being lowered, as expected.  The addition of inertia results in the 
membrane oscillating around the steady state, in a way that resembles
the case of the overdamped spring model (\ref{e:spring}).  
The inertia is responsible for the lowering of $\lambda^{*}_{w}$ with respect 
to $\lambda^{*}_{s}$, even in the limiting case of small aspect ratio $\epsilon = 0$, 
as reported by Flores \cite{flores}.  However, the aspect ratio has a stronger effect 
on the lowering of $\lambda^{*}_{w}$.

\begin{figure}
\centering
\includegraphics[scale=0.4,angle=270]{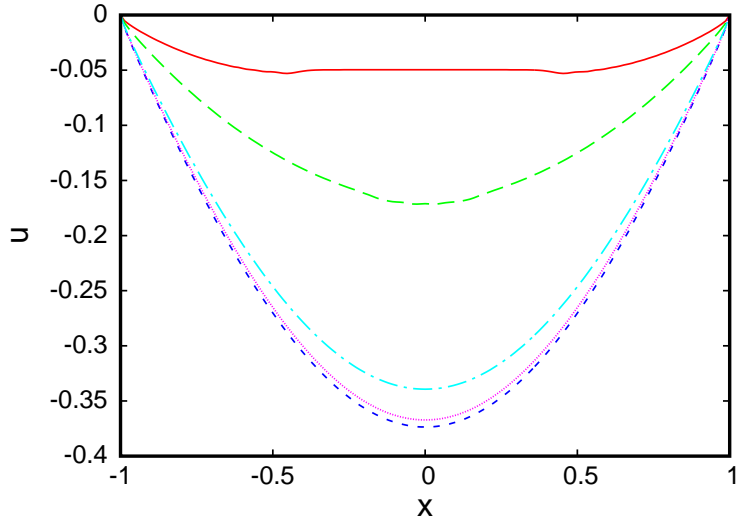}
\caption{(Color online)  Numerical solution of potential equation
(\ref{e:transpot}) and forced, damped wave equation (\ref{e:ueqntwave}) for 
$\lambda = 0.34$,
$\epsilon = 0.1$ and $\gamma = 0.7$.  $t=0.5$: red (solid) line, $t=1.0$: 
green (long dashed line), $t=3.0$: blue (short dashed) line, $t=5.0$: pink 
(dotted) line,
$t=10.0$: light blue (dot-dashed) line.}
\label{f:wave}
\end{figure}

The oscillatory approach of $u$ to the steady state when the displacement
$u$ is governed by the forced, damped wave equation (\ref{e:ueqntwave}) 
is illustrated in Figure \ref{f:wave}.  The parameter values $\lambda = 0.34$, 
$\epsilon = 0.1$ and $\gamma = 0.7$ were chosen so that the evolution is just 
below the critical $\lambda^{*}_{w} = 0.34468$.  The evolution is shown until 
the steady state is reached by $t=10$.  The profile reaches a maximum depth for 
$t \approx 3$ and then oscillates back up.  Before this time, the profiles are 
monotonically increasing in depth.  After $t = 3$ the profiles monotonically decrease 
in depth until the steady state is reached.  Thus, the behavior is similar to that of 
a heavily damped spring, as in the model (\ref{e:spring}) as reported by Rocha 
{\em et al} \cite{rocha}.

\begin{figure}
\centering
\includegraphics[scale=0.5,angle=270]{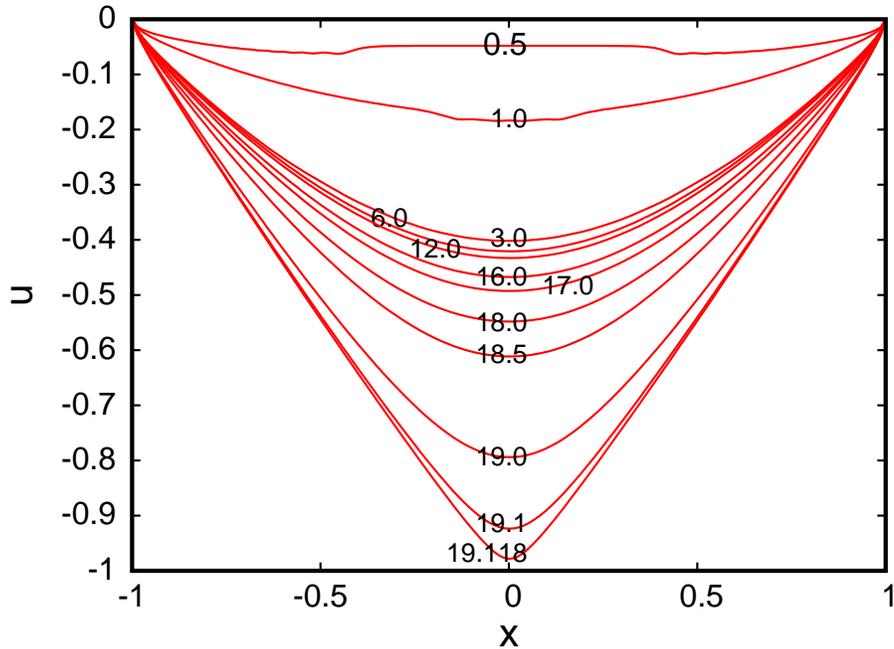}
\caption{(Color online)  Numerical solution of potential equation 
(\ref{e:transpot})
and forced, damped wave equation (\ref{e:ueqntwave}) for $\lambda = 0.327$, 
$\epsilon = 0.2$ and $\gamma = 0.7$.  The numbers on the curves are the time 
$t$ for the solution.}
\label{f:waveunstab}
\end{figure}

The contrasting evolution when quenching occurs is illustrated in Figure 
\ref{f:waveunstab}.  The parameter values were chosen just above the critical 
$\lambda^{*}_{w} = 0.3251$, with $\lambda = 0.327$, $\epsilon = 0.2$ and 
$\gamma = 0.7$. For these parameter values, $\lambda^{*}_{s}=0.32738$.
First, the profiles move on a fast time scale and approach the 
steady state corresponding to the static critical value.  There is then a slow 
motion away from that steady state until the depth has increased far enough 
that  the profiles can move on a fast time scale towards $u=-1$.  
This is similar to the observations of Rocha {\em et al} \cite{rocha} for the model 
(\ref{e:spring}).  Kavallaris {\em et al} \cite{kavallaris} obtained oscillations 
around the steady state and later approach to touchdown for values of $\lambda$
close to, but smaller than the critical static value. The oscillations 
are explained by the fact that their model corresponds to the regime in
which inertial forces dominate.
For our problem, in the third stage the displacement $u$ rapidly 
approaches quenching, at which point the numerical solution breaks down.

\section{Conclusions}

The static and dynamical behaviour of a flexible membrane driven by an 
electric field in a MEMS device has been investigated.  This evolution is governed by
a potential equation for the electric field with a nonlinear boundary condition
giving the membrane profile.  This moving boundary problem was transformed into
a boundary value problem on a fixed, rectangular domain, which was then 
investigated numerically due to the complexity of these equations.  
One of the findings is that the bifurcation curve has a single turning point 
with a shape which is qualitatively similar to that obtained in the limiting 
case of vanishing aspect ratio.  The dynamical evolution of the membrane was 
investigated by replacing the static membrane equation with both a forced heat equation 
and a forced, damped wave equation.  It was found that there is a critical value of the 
applied voltage for which the membrane does not settle to a steady state, but ``quenches,'' 
that is, it hits the bottom of the MEMS device, at which point the governing equations 
become invalid. In the case of the forced heat equation the dynamical critical 
value was found to be equal to the static critical value. In the case of the
forced, damped wave equation, the dynamical critical value
is lower than that for  the static problem. The numerical results show that
the dynamical and static critical values are different due to the inertial
forces. However, the aspect ratio is more important in the determination
of the dynamical critical value.  This is due to the membrane oscillating in its evolution.  
These results show the increased complexity which arises from more realistic models of the 
MEMS device.

\vskip.3cm

{\bf Acknowledgments.} The financial support of CONACyT through 
Proyecto de Grupo 133036-F is gratefully acknowledged.


\begin{thebibliography}{99}

\bibitem{desai} A.V. Desai, J.D. Tice, C.A. Apblett and P.J.A. Kenis,
``Design consideration for electrostatic microvalves with applications in
poly(dimethylsiloxane)-based microfluidics,''
{\em Lab Chip,} {\bf 12}, 1078--1088 (2012).

\bibitem{rocha} L.A. Rocha, E. Cretu and R.F. Wolffenbuttel,
``Behavioural analysis of the pull-in dynamical transition,''
{\em J. Microelectromech.\ Microeng.,} {\bf 14}, S37--S42 (2004).

\bibitem{zhang} W.-M. Zhang, H. Yan, Z.-K. Peng and G. Meng,
``Electrostatic pull-in instability in MEMS/NEMS: A review,''
{\em Sensors and Actuators A: Physical,} {\bf 214}, 187--218 (2014).
 
\bibitem{pelesko} J.A. Pelesko and D.H. Bernstein,
{\em Modelling MEMS and NEMS,} Chapman and Hall/CRC (2003).

\bibitem{flores2} G. Flores, G. Mercado, J.A. Pelesko and N.F. Smyth,
``Analysis of the dynamics and touchdown in a model of electrostatic MEMS,''
{\em SIAM J. Appl.\ Math.,} {\bf 67}, 434--446 (2007).

\bibitem{pelesko2} J.A. Pelesko and T.A. Driscoll,
``The effect of the small-aspect-ratio approximation on canonical
electrostatic MEMS models,'' {\em J. Engng.\ Math.,} {\bf 53}, 239-252 (2005).

\bibitem{brubaker} N.D. Brubaker, J.A. Pelesko, ``Nonlinear effects on
canonical MEMS models,'' {\em Euro. Jnl. of Applied \ Mathematics,}
{\bf 22}, 455--470, (2011).

\bibitem{lindsay} A.E. Lindsay and M.J. Ward, ``Aymptotics of some nonlinear
eigenvalue problems for a MEMS capacitor: Part I: Fold point asymptotics,''
{\em Meth. and Appl. of Analysis,} {\bf 15}, 297--326, 2008. 

\bibitem{wei} J. Wei and D. Ye, ``On MEMS equation with fringing field,''
{\em Proc.\ A.M.S.,} {\bf 138}, 1693--1699 (2010).

\bibitem{liu} Z. Liu and X. Wang,
``On a parabolic equation in MEMS with fringing field,''
{\em Arch.\ Math.,} {\bf 98}, 373--381 (2012).

\bibitem{flores} G. Flores,
``Dynamics of a damped wave equation arising from MEMS,''
{\em SIAM J. Appl.\ Math.,} {\bf 74}, 1025--1035 (2014).

\bibitem{chang} P.H. Chang and H.A. Levine, ``The quenching of solutions of
semilinear hyperbolic equations,'' {\em SIAM J. \ Math. Anal.,} {\bf 12},
893--903, 1982.

\bibitem{smith} R.A. Smith, ``A quenching problem in several dimensions,''
{\em SIAM. J. \ Math. Anal.,} {\bf 20} 1081--1094, 1989.

\bibitem{kavallaris} N.I. Kavallaris, A.A. Lacey, C.V. Nikolopoulos and
D.E. Tzanetis, ``A hyperbolic nonlocal problem modelling MEMS technology,''
{\em Rocky Mountain J. \ Math.,} {\bf 41}, 505--534, 2011.

\bibitem{siddique} J.I. Siddique, R. Deaton, E. Sabo and J.A. Pelesko,
``An experimental investigation of the theory of electrostatic deflections,''
{\em J. Electrostatics}, {\bf 69}, 1--6 (2011).

\bibitem{lauren} P. Lauren\c{c}ot and C. Walker,
``A stationary free boundary problem modeling electrostatic MEMS,''
{\em Arch.\ Rational Mech.\ Anal.,} {\bf 207}, 139--158 (2013).

\bibitem{escher} J. Escher, P. Lauren\c{c}ot and  C. Walker,
``A parabolic free boundary problem modeling electrostatic MEMS,''
{\em Arch.\ Rat.\ Mech.\ Anal.,} {\bf 211}, 389--417 (2014).






\end{thebibliography}
\end{document}